\documentclass[11pt]{amsart}
\usepackage{fullpage}
\usepackage {amsmath}
\usepackage {amsthm}
\usepackage {enumerate}
\usepackage{amssymb}
\usepackage{amsfonts}
\usepackage{anysize}
\usepackage{setspace}

\marginsize{3.15cm}{3.15cm}{1.6cm}{3cm}

\newcommand{\N}{\mathbb N}

\newcommand{\R}{\mathbb R}

\newcommand{\A}{\mathcal{A}}

\numberwithin{equation}{section}

\theoremstyle{plain}

\newtheorem{theorem}[equation]{Theorem}

\newtheorem{corollary}[equation]{Corollary}
\newtheorem{lemma}[equation]{Lemma}
\newtheorem{question}[equation]{Question}
\newtheorem{problem}[equation]{Problem}

\theoremstyle{definition}

\theoremstyle{definition}
\newtheorem{remark}[equation]{Remark}

\theoremstyle{definition}

\theoremstyle{definition}
\newtheorem{example}[equation]{Example}

\theoremstyle{definition}

\author{{\bfseries S. Garc\'ia-Ferreira}}
\address{Centro de Ciencias Matem\'aticas\\
         Universidad Nacional Aut\'onoma de M\'exico\\
				 Campus Morelia\\
         Apartado Postal 61-3, Santa Mar\'ia, 58089, Morelia, Michoac\'an, M\'exico.}
\email{sgarcia@matmor.unam.mx}

\title[Categorical properties on the  hyperspace of convergent sequences]{\scshape\bfseries Categorical properties on the  hyperspace of nontrivial convergent sequences}

\author{{\bfseries R. Rojas-Hern\'andez}}
\email{satzchen@yahoo.com.mx}

\author{ Y. F. Ortiz-Castillo }
\address{Instituto de Matem\'atica e Estat\'istica, Universidade de S\~ao Paulo \\ Rua do Mat\~ao, 1010, CEP 05508-090, S\~ao Paulo, Brazil}
\email{jazzerfoc@gmail.com}



\date{}

\thanks{Research of the first-named author was supported
by  CONACYT grant no. 176202 and PAPIIT grant no. IN-101911. The research of the third listed author  was support from CNPq (Brazil) - ``Bolsa de Produtividade em Pesquisa, processo 307130/2013-4'', CNPq Universal 483734/2013-6. }

\begin{document}

\begin{abstract}  In this paper, we shall study categorial properties of the
hyperspace of all nontrivial convergent sequences $\mathcal{S}_c(X)$ of a Fre\'ech-Urysohn space  $X$, this hyperspace  is equipped with the Vietoris topology.
We mainly  prove that  $\mathcal{S}_c(X)$ is meager whenever $X$ is a crowded space, as a corollary we obtain that if $\mathcal{S}_c(X)$ is Baire, the $X$ has a dense subset of isolated points.
As an interesting example $\mathcal{S}_c(\omega_1)$ has the Baire property, where $\omega_1$ carries the order topology (this  answers a question from \cite{sal-yas}). We can give more examples like this one by proving that the  Alexandroff duplicated $\mathcal{A}(Z)$ of a space $Z$ satisfies that $\mathcal{S}_c(\mathcal{A}(Z))$ has the Baire property, whenever    $Z$ is a $\Sigma$-product of completely metrizable spaces and $Z$ is crowded. Also we show that if
$\mathcal{S}_c(X)$ is pseudocompact, then  $X$ has a relatively countably compact dense subset of isolated points, every finite power of $X$ is pseudocompact, and every $G_\delta$-point in $X$ must be isolated. We also establish several topological properties  of the hyperspace of nontrivial convergent sequences of countable  Fre\'ech-Urysohn spaces with only one non-isolated point.
\end{abstract}

\maketitle

\section{Introduction}

 All our spaces will be Tychonoff (completely regular and Hausdorff). The letters  $\mathbb{P}$ and  $\N$ will denote the irrational numbers and  the natural numbers, respectively. The positive natural numbers will be denoted by  $\N^+$. The Greek letter $\omega$ stands for the first infinite cardinal number, $\omega_1$ stands for the first uncountable cardinal number endowed with the order topology. If $A, B \subseteq \mathbb{N}$, then $A \subseteq^* B$ means that $A \setminus B$ is finite.

\medskip

For a topological space $X$, $\mathcal{CL}(X)$ will denote the set of all nonempty closed subsets of $X$. For a nonempty family $\mathcal{U}$ of subsets of $X$ let
\begin{center}
$\left\langle \mathcal{U} \right\rangle = \{F \in \mathcal{CL}(X) : F \subseteq \bigcup{\mathcal{U}} \textnormal{ and } F \cap U \not= \emptyset\}.$
\end{center}
If $\mathcal{U} = \{U_1,\ldots,U_n\}$, in some convenient cases,  $\left\langle \mathcal{U} \right\rangle$ will be denoted by  $\left\langle U_1,\ldots,U_n \right\rangle$. A base for the \textit{Vietoris topology} on $\mathcal{CL}(X)$ is the family of all sets of the form $\left\langle \mathcal{U}\right\rangle$, where $\mathcal{U}$ runs over all nonempty finite families of nonempty open subsets of $X$. In the sequel, any subset $\mathcal{D} \subseteq \mathcal{CL}(X)$ will carry the relative Vietoris topology as a subspace of $\mathcal{CL}(X)$. Given $\mathcal{D} \subseteq \mathcal{CL}(X)$ and a nonempty family $\mathcal{U}$ of subsets of $X$ we let $\left\langle \mathcal{U}\right\rangle_\mathcal{D} = \left\langle \mathcal{U}\right\rangle \cap \mathcal{D}$. For simplicity, if there is no possibility of confusion we simply write  $\left\langle \mathcal{U}\right\rangle$ instead of $\left\langle \mathcal{U}\right\rangle_\mathcal{D}$.
All topological notions whose definition is not included here should be understood as in \cite{eng}.

\medskip

Some of the most studied hyperspaces on a space $X$ have been
$$
\mathcal{K}(X) = \{K \in \mathcal{CL}(X) : K \textnormal{ is compact}\} \textnormal{ and } \mathcal{F}(X) = \{F \subseteq X : F \textnormal{ is finite and} \ F \neq \emptyset \},
$$
see for instance the survey paper \cite{gu}.
In the article \cite{jvpp}, the authors considered the hyperspace consisting of all finite subsets together with all  the Cauchy sequences  without limit point of a metric space.
In a different context,  the authors of the paper \cite{peng} consider the set $\mathcal{F}_S(X)$ of all convergent sequences of a space $X$, and study the existence of a metric $d$ on the set $X$ such that $d$ metrizes all subspaces of $X$ which belong to $\mathcal{F}_S(X)$; that is, the restriction of $d$ to $A$ generates the subspace topology on $A$ for every $A \in \mathcal{F}_S(X)$ (these kind of problems have been analyzed in \cite{arh}).

\medskip

The hyperspace of nontrivial convergent sequences was studied in  \cite{sal-yas} and, recently, in  \cite{may-pat-rob} and \cite{salrey}. This hyperspace is formally defined as follows:

\smallskip

 Given a space $X$, a \textit{nontrivial convergent sequence} of $X$ is a subset $S \subseteq X$ such that  $\left|S\right| = \omega$, $S$ has a unique non-isolated point, which  will be denote by  $x_S$, and $\left|S \setminus U\right| < \omega$ for each neighborhood $U$ of $x_S$. By using  this notion, we define
$$\mathcal{S}_c(X) = \{S \in \mathcal{CL}(X) : S \textnormal{ is a nontrivial convergent sequence}\}.$$
It was pointed out in \cite{may-pat-rob} that the family of all subsets of $\mathcal{S}_c(X)$ of the form $\left\langle \mathcal{U}\right\rangle$, where $\mathcal{U}$ is a finite family of pairwise disjoint open subsets of $X$,  is a base for the Vietoris topology on $\mathcal{S}_c(X)$. We will refer to  this family as the {\it canonical basis} of $\mathcal{S}_c(X)$ and  its elements will be named as {\it canonical open sets}.

\medskip

Throughout this paper, it is evident why we shall only consider  non-discrete  Fr\'echet-Urysohn spaces. Thus, the hyperspace $\mathcal{S}_c(X)$ will be always nonempty.

\medskip

A fundamental task in the study of the hyperspace $\mathcal{S}_c(X)$ is to determine its topological relationship with  the base space $X$ and vise versa.
For instance, the connection  between the connectedness in $X$ and $\mathcal{S}_c(X)$  was explore  in \cite{sal-yas}, \cite{may-pat-rob} and \cite{salrey}.
The category property on $\mathcal{S}_c(X)$ was another of the topological properties studied in \cite{sal-yas}, it was proved in there  that $\mathcal{S}_c(X)$ is never a Baire space when the space $X$ is crowded, and that $\mathcal{S}_c(X)$ is even meager if, in addition, $X$ is second countable. Following this direction,  our main purpose of this article is to continue studying the category properties on $\mathcal{S}_c(X)$.

\medskip

The first section is devoted to study the hyperspace of countable Fre\'ech-Urysohn spaces with just one accumulation point. We show that if $X$ is such a space, then  $\mathcal{S}_c(X)$ is homeomorphic to the rational numbers iff $X$ has a countable base.

 The last section is devoted to study categorical  properties of the hyperspace of nontrivial convergent sequences.  Concerning this topic,  the authors of \cite[Questions 3.4 and 3.5 ]{sal-yas}  asked about the meager property in $\mathcal{S}_c(X)$ when $X$ is a metrizable crowded space, and about the meager property of the space $\mathcal{S}_c(\omega_1)$. We answer both questions by showing that $\mathcal{S}_c(X)$ is meager whenever $X$ is crowded, and that $\mathcal{S}_c(\omega_1)$ has the Baire property. Finally, we prove that if
$\mathcal{S}_c(X)$ is pseudocompact, then  $X$ has a relatively countably compact dense subset of isolated points, every finite power of $X$ is pseudocompact, and every $G_\delta$-point in $X$ must be isolated. The last result of this paper is that  if $Z$ is a $\Sigma$-product of completely metrizable spaces and $Z$ is crowded, then $\mathcal{S}_c(\mathcal{A}(Z))$ has the Baire property, where $\mathcal{A}(Z)$ is the Alexandroff duplicated of $Z$. For further research  we list some open questions related to the results of the paper.

\section{Countable Fr\'echet-Urysohn spaces}

The main result of this section is to give  two non-homeomorphic spaces $X$ and $Y$ such that $\mathcal{S}_c(X)$ and $\mathcal{S}_c(Y)$ are homeomorphic.

\medskip

Since $\mathcal{S}_c(X)$ is, in general, a dense proper subset of  $\mathcal{CL}(X)$, it cannot be compact. However, we can say more:

\begin{lemma}\label{nolocallycom} For very space $X$, $\mathcal{S}_c(X)$ is nowhere locally compact.
\end{lemma}

\proof  Fix $S \in \mathcal{S}_c(X)$ and let $O$ be a neighborhood of $S$. Choose   a canonical open set $\left\langle U_1,...,U_n  \right\rangle$ so that $S \in \left\langle cl_X(U_1),...,cl_X(U_n)  \right\rangle \subseteq O$,  the sets  $\{cl_X(U_1),...,cl_X(U_n)\}$ are pairwise disjoint, $\{x_i\} = S \cap cl_X(U_i)$, for each $1 \leq i < n$, and  $x_S \in U_n$. Choose a local base $\mathcal{B}$ at $x_S$ and one more $\mathcal{B}_i$ at $x_i$, for each $1\leq i < n$. Notice that $\left\langle cl_X(U_1 \cap B_1),...,cl_X(U_{n-1} \cap B_{n-1}),cl_X(U_n \cap B)) \right\rangle$ is a closed subsets of $\mathcal{S}_c(X)$ for each $(B_1,...,B_{n-1},B) \in \mathcal{B}_1 \times \mathcal{B}_2\times \cdots \mathcal{B}_{n-1} \times \mathcal{B}$. So, $$
  \{\left\langle cl_X(U_1 \cap B_1),...,cl_X(U_{n-1} \cap B_{n-1}),cl_X(U_n) \cap B) \right\rangle :
  $$
  $$
  (B_1,...,B_{n-1},B) \in \mathcal{B}_1 \times \mathcal{B}_2\times \cdots \mathcal{B}_{n-1} \times \mathcal{B} \}
   $$
   is a family of closed subsets of $\left\langle \mathcal{U}\right\rangle$ with the finite intersection property. But
   $$
   \bigcap\{\left\langle cl_X(U_1 \cap B_1),...,cl_X(U_{n-1} \cap B_{n-1}),cl_X(U_n) \cap B) \right\rangle :
   $$
   $$
   (B_1,...,B_{n-1},B) \in \mathcal{B}_1 \times \mathcal{B}_2\times \cdots \mathcal{B}_{n-1} \times \mathcal{B}  \} = \emptyset.
   $$
   Therefore, neither  $\left\langle cl_X(U_1),...,cl_X(U_n)  \right\rangle$ nor $cl_X(O)$ can be  compact.
   \endproof

The spaces that we are looking for will be countable with only one non-isolated point. To deal with this kind of spaces we shall need the following terminology:

\smallskip

To each free filter $\mathcal{F}$ on $\N$ we associate the space $\xi(\mathcal{F})$ whose underlying set is $\N \cup \{\mathcal{F}\}$, all elements of $\N$ are isolated and the neighborhoods of $\mathcal{F}$ are of the form $A \cup \{\mathcal{F}\}$ where $A \in \mathcal{F}$.  It is evident that the space $\xi(\mathcal{F})$ is zero-dimensional and Hausdorff for every free filter  $\mathcal{F}$ on $\N$. A free filter  $\mathcal{F}$ on $\N$ is said to be a {\it Fr\'echet-Urysohn filter} if the space $\xi(\mathcal{F})$ is Fr\'echet-Urysohn (there are plenty of this kind of filters, see \cite{guz}). The most simplest example of a Fr\'echet-Urysohn filter is the {\it Fr\'echet filter} $\mathcal{F}_r = \{ A \subseteq \N : |\N \setminus A| < \omega  \}$, and notice that $\xi(\mathcal{F}_r)$ is a convergent sequence with its limit point. One more example of a  Fr\'echet-Urysohn filter with a countable base is the filter $\mathcal{P} =  \{ A \subseteq \mathbb{N} : \exists m \in \mathbb{N} \forall n \geq m ( P_n \subseteq A) \}$, where $\{ P_n : n \in \N \}$ is a partition of $\N$ in infinite subsets. At this point, we can say that these two Fr\'echet-Urysohn filters are the only ones, up  to homeomorphism, with a countable base. These two filters can be characterized as follows.

\begin{theorem}\label{SCFCB} Let $\mathcal{F}$ be a Fr\'echet-Urysohn filter. Then
 $\mathcal{S}_c(\xi(\mathcal{F}))$ is homeomorphic to $\mathbb{P}$ iff $\mathcal{F}$ has a countable base.
\end{theorem}

\proof Necessity.
If $\mathcal{S}_c(\xi(\mathcal{F}))$ is homeomorphic to $\mathbb{P}$, then $\mathcal{S}_c(\xi(\mathcal{F}))$ is second countable. Hence, by applying \cite[Theorem 6.5]{may-pat-rob}, we obtain that $\xi(\mathcal{F})$ is also second countable. Thus, we conclude that the filter $\mathcal{F}$ has a countable base.

\medskip

Sufficiency. Assume  that $\mathcal{B}$ is a countable base of $\mathcal{F}$. It then follows  from  \cite[Corollary 6.15]{may-pat-rob} that the  space $\mathcal{S}_c(\xi(\mathcal{F}))$ is separable.  Since every open subset of $\xi(\mathcal{F})$ is closed, we may apply \cite[Prop. 3.2]{may-pat-rob} to see that the canonical basis of $\mathcal{S}_c(\xi(\mathcal{F}))$ consists of clopen subsets. So,
we also obtain that $\mathcal{S}_c(\xi(\mathcal{F}))$ is zero-dimensional. In order to see that $\mathcal{S}_c(\xi(\mathcal{F}))$ is homeomorphic to $\mathbb{P}$, in virtue of \cite[Theorem 1.9.8]{jvm}, it is enough to show that it is nowhere locally compact and completely metrizable. Indeed,
by Lemma \ref{nolocallycom}, we know that $\mathcal{S}_c(\xi(\mathcal{F}))$ is nowhere locally compact.
Let us show that $\mathcal{S}_c(\xi(\mathcal{F}))$ is completely metrizable. Enumerate $\mathcal{B}$ as $\{B_n : n \in \N\}$ and, without loss of generality, assume that $B_{n+1} \subseteq B_n$ for each $n \in \N$. Let $E_n = B_n \setminus B_{n+1}$ for each $n \in \N$. Without loss of generality, we may assume that $\{ E_n : n \in \N\}$ is a partition of $\N$. Consider the map $f : \xi(\mathcal{F}) \to \R$ defined as $f(m) = 1/2^n$ for $m \in E_n$ and $f(\mathcal{F}) = 0$. Define a metric $d$ on $\xi(\mathcal{F})$ as follows: For $x,y \in \xi(\mathcal{F})$, if $x,y \in E_n$ for some $n \in \N$ and $x \not= y$, then we let $d(x,y) = 1/2^n$ and we define $d(x,y) = \left|f(x) - f(y)\right|$ otherwise. Notice that $d$ is a complete metric compatible with the topology of $\xi(\mathcal{F})$. It is well-known that the Hausdorff metric induced by $d$ on $\mathcal{CL}(\xi(\mathcal{F}))$ is also complete (see \cite{smith}). Thus, we obtain that  $\mathcal{CL}(\xi(\mathcal{F}))$ is completely metrizable. To show that $\mathcal{S}_c(\xi(\mathcal{F}))$ is completely metrizable, it is enough to prove that $\mathcal{S}_c(\xi(\mathcal{F}))$ is a $G_\delta$-set in $\mathcal{CL}(\xi(\mathcal{F}))$. For each $B \in \mathcal{B}$ we let $O_B = \{S \in \mathcal{CL}(\xi(\mathcal{F})) : \left|S \setminus B\right| < \omega\}$ and note that this set is open in $\mathcal{CL}(\xi(\mathcal{F}))$; indeed,
\begin{center}
$O_B = \bigcup\{\left\langle \{B\} \cup \{\{x\} : x \in F\} \right\rangle : F \in \mathcal{F}(\xi(\mathcal{F}) \setminus B) \cup \{\emptyset\}\} \cup \mathcal{F}(\xi(\mathcal{F}) \setminus B).$
\end{center}
On the other hand, since $\mathcal{F}(\xi(\mathcal{F}))$ is $F_\sigma$ in $\mathcal{CL}(\xi(\mathcal{F}))$, it follows that $\mathcal{S}_c(\xi(\mathcal{F})) = \bigcap\{O_B : B \in \mathcal{B}\} \setminus \mathcal{F}(\xi(\mathcal{F}))$ is a $G_\delta$-subset of $\mathcal{S}_c(\xi(\mathcal{F}))$.
\endproof

\begin{example}\label{NHSC}
There are two spaces $X$ and $Y$ such that $\mathcal{S}_c(X)$ is homeomorphic to $\mathcal{S}_c(Y)$ but $X$ is not homeomorphic to $Y$.
\end{example}

\proof
By  Theorem \ref{SCFCB}, we know that $\mathcal{S}_c(\xi(\mathcal{F}_r))$ is homeomorphic to $\mathcal{S}_c(\xi(\mathcal{P}))$, but it is clear that  the spaces $\xi(\mathcal{F}_r)$ and $\xi(\mathcal{P})$ cannot be homeomorphic.
\endproof

Both spaces $X$ and $Y$ considered in the previous example have a dense set isolated points. However, we still do not know any counterexample in the realm of crowded (Fr\'echet-Urysohn) spaces.

\begin{problem}\label{qfu}
Find two non-homeomorphic crowded spaces $X$ and $Y$ such that  $\mathcal{S}_c(X)$ and $\mathcal{S}_c(Y)$ are homeomorphic.
\end{problem}

Addressing to Question \ref{qfu}, we would like to make some comments about  the following class of sequential,  non-Fr\'echet-Urysohn crowded spaces:

\medskip

The symbol $FF(\N)$ will denote the set of all free filters on $\N$ and $Seq =  \bigcup_{n  \in N}\N^n$.
If $s \in Seq$ and $n \in \N$, then the {\it concatenation} of $s$ and $n$ is the function $s^\frown n = s \cup \{(dom(s),n)\}$. For a function $\delta: Seq \to FF(\N)$, we define a topology $\tau_\delta$ on $Seq$ by defining $V \in \tau_\delta$ iff $\{ n  \in \N : s^\frown n \in V \} \in \delta(s)$ for every $s \in V$. It is well-known that  $ Seq(\delta) = (Seq,\tau_\delta)$ is a extremally disconnected, zero dimensional Hausdorff space for every function $\delta: Seq \to FF(\N)$. Besides, $Seq(\delta)$ is a  sequential spaces provided that  $\delta(s)$ is a Fr\'echet-Urysohn filter for each $s \in Seq$. It is not hard to see that
$Seq(\delta)$ cannot be Fr\'echet-Urysohn.

\begin{question}\label{qfug}
Are   $\mathcal{S}_c(Seq(\mathcal{F}_r))$ and $\mathcal{S}_c(Seq(\mathcal{P}))$ homeomorphic?
\end{question}

After having Theorem \ref{SCFCB}, we shall give next a necessary conditions when the space $\mathcal{S}_c(Seq(\mathcal{F}))$ is Baire.

\medskip

Following A.V. Arhangel'skii \cite{arh81}, we say that $x \in X$ is an $\alpha_2$-{\it point} if for every family $\{\ S_n : n \in \mathbb{N} \}$ of sequences converging to $x$ there is $S \in  \mathcal{S}_c(X)$ converging to $x$ such that $|S \cap S_n| = \omega$ for all $n \in \mathbb{N}$. We shall say that  a filter $\mathcal{F}$ is an $\alpha_2$-{\it filter} if every point of the space $\xi(\mathcal{F})$ is an $\alpha_2$-point.

\begin{lemma}\label{secondcategory} Let $\mathcal{F}$ be a Fr\'echet-Urysohn filter,  $S \in \mathcal{S}_c(\xi(\mathcal{F}))$ and $n \in \mathbb{N}$. Then $\mathcal{D}_S^n = \{ T \in \mathcal{S}_c(\xi(\mathcal{F})) : |T \cap S| \geq n\}$  is a dense open subset of $ \mathcal{S}_c(\xi(\mathcal{F}))$.
\end{lemma}

\proof Fix $S \in \mathcal{S}_c(\xi(\mathcal{F}))$ and $n \in \mathbb{N}$. Assume that $\langle  \mathcal{U} \cup \{A\} \rangle$ is a canonical open set of $\mathcal{S}_c(\xi(\mathcal{F}))$, where  $\mathcal{U}$ consists of singletons which are elements of $\mathbb{N}$ and $A \in \mathcal{F}$.
Since $S \subseteq^* A$, then we can easily find $T \in \langle  \mathcal{U} \cup \{A\} \rangle$ so that $|T \cap S| \geq n$. To prove that $\mathcal{D}_S^n$ is open fix $T \in \mathcal{D}_S^n$. Then choose a family of singletons  $\mathcal{U}$ consisting of  $n$ elements of $T \cap S$. It is then clear that
$T \in \langle  \mathcal{U} \cup \{\mathbb{N}\} \rangle \subseteq \mathcal{D}_S^n$.
\endproof

\begin{theorem}\label{alpha} Let $\mathcal{F}$ be a Fr\'echet-Urysohn filter.  If the hyperspace
 $\mathcal{S}_c(\xi(\mathcal{F}))$ is not meager, then $\mathcal{F}$ is an $\alpha_2$-filter.
\end{theorem}

\proof Let $\{ S_m : m \in \mathbb{N} \}$ be a countable subset of $\mathcal{S}_c(\xi(\mathcal{F}))$.
 According to Lemma \ref{secondcategory},  the set  $\mathcal{D}_{S_m}^n$ is open and dense for each $n, m \in \mathbb{N}$. Hence, we can find
 $T \in \bigcap_{n, m \in \mathbb{N}}\mathcal{D}_{S_m}^n$. It is then clear that $|T \cap S_m| = \omega$ for every $m \in \mathbb{N}$.
\endproof

To give an other example of a Frech\'et-Urysohn filter we need
the characterization of the Frech\'et-Urysohn filters given in
\cite{si98}:
 A free filter $\mathcal{F}$ is Frech\'et-Urysohn iff there
is an $AD$-family $\mathcal{A}$ maximal with respect to the
following properties:
\begin{enumerate}
\item $F \in \mathcal{F}$ iff $A \subseteq^* F$ for all $A \in
\mathcal{A}$, and

\item $\mathcal{A}$ is an $AD$-family.
\end{enumerate}
As a consequence, we have that if $S \in \mathcal{S}_c(\xi(\mathcal{F}))$, then $|S \cap A| = \omega$ for some $A \in \mathcal{A}$. When $\mathcal{F}$ is a Frech\'et-Urysohn and $\mathcal{A}$
is the $AD$-family given by its characterization, we shall put $\mathcal{F} = \mathcal{F}_{\mathcal{A}}$. If $\mathcal{Q}$ is a partition of $\N$, then $\mathcal{F}_{\mathcal{Q}}$ is the well-known $FAN$-{\it filter}.
Since the $FAN$-filter  $\mathcal{F}_{\mathcal{Q}}$ cannot be an $\alpha_2$-filter, it follows from the previous theorem that $\mathcal{S}_c(\xi(\mathcal{F}_{\mathcal{Q}}))$ is meager.

\medskip

The following problem is then natural.

\begin{problem}
Determine the Fr\'echet-Urysohn filters $\mathcal{F}$ on $\N$ for which the space $\mathcal{S}_c(\xi(\mathcal{F}))$ is Baire.
\end{problem}

\section{Category in $\mathcal{S}_c(X)$}

Theorem \ref{SCFCB} provides an example of a space $X$ for which $\mathcal{S}_c(X)$ is Baire.  In what follows, we shall describe more examples of spaces for which  $\mathcal{S}_c(X)$ has this property. To have this done we list some easy facts and introduce some useful notation.

\begin{remark}\label{ROD}
Let Y be a non-discrete subspace of $X$, then we have that:
\begin{enumerate}
\item[$(a)$] $Y$ is dense in $X$ iff $\mathcal{S}_c(Y)$ is dense in $\mathcal{S}_c(X)$.

\item[$(b)$] $Y$ is open in $X$ iff $\mathcal{S}_c(Y)$ is open in $\mathcal{S}_c(X)$.
\end{enumerate}
As a consequence, if $Y$ is open and dense in $X$, then $\mathcal{S}_c(Y)$ is Baire iff $\mathcal{S}_c(X)$ is Baire.

\smallskip

 For a space $X$, $\mathcal{S} \subseteq \mathcal{S}_c(X)$ and $D \subseteq X$, we define
  $$
  \mathcal{G}(\mathcal{S},D) = \{S \cup F : S \in \mathcal{S} \textnormal{ and } F \in \mathcal{F}(D) \cup \{\emptyset\}\}.
  $$
The set $\mathcal{G}(\mathcal{S},X)$ will be simply denoted by $\mathcal{G}(\mathcal{S})$. If  $\mathcal{S} \subseteq \mathcal{S}_c(X)$, then we have that $\mathcal{S} \subseteq \mathcal{G}(\mathcal{S},D)$ for all $D \subseteq X$.  The following properties can be easily verified.
\begin{enumerate}
\item[$(c)$] If $\mathcal{S} \subseteq \mathcal{S}_c(X)$ is dense, then  $\mathcal{G}(\mathcal{S},D)$ is dense in $ \mathcal{S}_c(X)$ for all $D \subseteq X$.

\item[$(d)$] If $Y \subseteq X$ is open and  $D \subseteq X$ is discrete, then  $\mathcal{G}(\mathcal{S}_c(Y),D)$ is open in  $ \mathcal{S}_c(X)$.
\end{enumerate}
\end{remark}

In the next theorem, we state several conditions on a space $X$ that guarantee the Baire property of the hyperspace $\mathcal{S}_c(X)$. We shall need the following lemmas.

\begin{lemma}\label{opendis} Let $X$ be a space such that its set $D$ of isolated points  is dense in $X$. If $Y$ is a nonempty  open subset of $X$ such that $\mathcal{S}_c(Y)$ is Baire, then
 $\mathcal{G}(\mathcal{S}_c(Y),D)$ is Baire.
\end{lemma}

\proof
Assume that $\{\mathcal{D}_n : n \in \N\}$ is a family of open and dense subsets of $\mathcal{G}(\mathcal{S}_c(Y),D)$.
Let $\left\langle\mathcal{U}\right\rangle$ be a nonempty canonical open subset of $\mathcal{G}(\mathcal{S}_c(Y),D)$. Without loss of generality suppose that $\mathcal{U} = \mathcal{U}' \cup \mathcal{U}_D$ where $\mathcal{U}'$ is a family of  pairwise disjoint infinite open subsets of $Y$ and $\mathcal{U}_D$ is a family of singletons sets with  elements in $D$. Set $F = \bigcup \mathcal{U}_D$ and observe that $(S \cap Y) \cup F  = S$ for each  $S \in  \left\langle \mathcal{U}\right\rangle$. For each $n \in \N$, define $\mathcal{E}_n  = \{S \cap Y : S \in \mathcal{D}_n \cap \left\langle \mathcal{U}\right\rangle \}$.
It is straightforward to verify that  $\mathcal{E}_n$ is open and dense in $\left\langle \mathcal{U}'\right\rangle$,  for all  $n \in \N$. Since $\mathcal{S}_c(Y)$ is Baire and $\left\langle \mathcal{U}'\right\rangle$ is a nonempty open subset of $\mathcal{S}_c(Y)$, we can find $T \in \left\langle\mathcal{U}'\right\rangle \cap \big(\bigcap\{\mathcal{E}_n : n \in \N\}\big)$.  Then we have that  $T \cup F \in \left\langle\mathcal{U}\right\rangle \cap \big(\bigcap\{\mathcal{E}_n : n \in \N\}\big) \subseteq \left\langle\mathcal{U}\right\rangle \cap \big(\bigcap\{\mathcal{D}_n : n \in \N\}\big)$. Since $\left\langle \mathcal{U}\right\rangle$ was taken arbitrary, the space $\mathcal{G}(\mathcal{S}_c(Y),D)$ is Baire.
\endproof

We omit the proof of the following   well-known results.

\begin{lemma}\label{unionopenbaire} Let $X$ be a space.
\begin{enumerate}
\item If $\{ V_i : i \in I\}$ is a family of open Baire nonempty open subsets of $X$, then $\bigcup_{i \in I}V_i$ is also Baire.

\item If $X$ has  a family $\mathcal{U}$ consisting of pairwise disjoint   nonempty open meager subsets  and $\bigcup \mathcal{U}$ is dense in $X$, then $X$ is also meager.
\end{enumerate}
\end{lemma}

\begin{theorem}\label{CSCB} Let $X$ be a space such that its set $D$ of isolated points  is dense in $X$, and let  $\{X_\gamma: \gamma \in \Gamma\}$
be a family of  clopen subspaces of $X$. If the following conditions are satisfied:
\begin{enumerate}[i)]
\item the set $\{ x_S : S \in \bigcup_{\gamma \in \Gamma}\mathcal{S}_c(X_\gamma)\}$ is dense in $X \setminus D$,

\item $\mathcal{S}_c(X_\gamma)$ is Baire for each $\gamma \in \Gamma$,

\item	the family $\{\mathcal{G}(\mathcal{S}_c(X_\gamma),D) : \gamma \in \Gamma\}$ is pairwise disjoint, and

\item  $\mathcal{S}_c(X) = \bigcup\{\mathcal{G}(\mathcal{S}_c(X_\gamma)) : \gamma \in \Gamma\}$,
\end{enumerate}
 then $\mathcal{S}_c(X)$ is a Baire space.
\end{theorem}

\proof
First, we shall prove that $\mathcal{D} = \bigcup\{\mathcal{G}(\mathcal{S}_c(X_\gamma),D) : \gamma \in \Gamma\}$ is dense in $\mathcal{S}_c(X)$. Indeed,
 let $O = \left\langle\mathcal{U}\right\rangle$ be a nonempty canonical open subset of $\mathcal{S}_c(X)$. By condition $i)$, we may choose $S \in O$ so that $x_S \in X_\gamma$ for some   $\gamma \in \Gamma$. Notice  that $S \cap X_\gamma \in \mathcal{S}_c(X_\gamma)$. For each $U \in \mathcal{U}$ select a point $d_U \in D \cap U$, and consider the convergent sequence  $S_0 = (S \cap X_\gamma) \cup \{d_U : U \in \mathcal{U}\}$. Then we have that   $S_0 \in O \cap \mathcal{G}(\mathcal{S}_c(X_\gamma),D)$ and so $\mathcal{D} \cap O \not= \emptyset$. By clause $(d)$ of Remark \ref{ROD}, we obtain that $\mathcal{D}$ is open in $\mathcal{S}_c(X)$ and $\mathcal{D} = \bigcup\{\mathcal{G}(\mathcal{S}_c(X_\gamma),D) : \gamma \in \Gamma\}$.
By condition $ii)$ and Lemma \ref{opendis},  we know that $\mathcal{G}(\mathcal{S}_c(X_\gamma),D)$ is Baire for each $\gamma \in \Gamma$.  Thus, clause $(1)$ of Lemma \ref{unionopenbaire} implies that
 $\mathcal{D}$ is Baire and so $\mathcal{S}_c(X)$ is also Baire.
\endproof

\begin{example}\label{ex1}
 The space $\mathcal{S}_c(\omega_1)$ is Baire.
\end{example}

\proof
Let $D$ be the set of all isolated points of $\omega_1$ and
$$
Y = \{\alpha < \omega_1 : \alpha \textnormal{ is a limit ordinal and } (\beta, \alpha) \subseteq D \textnormal{ for some } \beta < \alpha\}.
$$
Set $X = Y \cup D$. Since $X$ is open and dense in $\omega_1$, according to Remark \ref{ROD} $(a)-(b)$, it is enough to show that $\mathcal{S}_c(X)$ is Baire. For each  $\alpha \in Y$ pick $\beta_\alpha$ so that $(\beta_\alpha, \alpha) \subseteq D$ and let $X_\alpha = (\beta_\alpha, \alpha]$. According to Theorem \ref{SCFCB}, we have that the space $\mathcal{S}_c(X_\alpha)$ is Baire for each $\alpha \in Y$. Thus, $D$ and the family $\{X_\alpha : \alpha \in Y\}$ satisfy all the conditions of Lemma \ref{CSCB}. Therefore,   $\mathcal{S}_c(X)$ is a Baire space.
\endproof

For the  description of  our next example, we recall that an {\it almost disjoint family}  $\mathcal{A}$ is an infinite family of infinite subsets of $\N$ such that $|A \cap B| < \omega$ for distinct $A, B \in \mathcal{A}$.
 The $\Psi$-{\it space} associated to an $AD$-family $\mathcal{A}$,  denoted by  $\Psi(\mathcal{A})$, is the space whose underlying set is $\N \cup \{ \{A\} : A \in \mathcal{A}\}$, $\N$ is discrete and, for each $A \in \mathcal{A}$,  $(A \setminus F) \cup \{A\}$ with $F \in \mathcal{F}(\N)$ is a basic neighborhood of $A$. It is easy to see that the space $\Psi(\mathcal{A})$ is always first countable and zero dimensional.

\begin{example}\label{ex2} For every  almost disjoint family $\mathcal{A}$ the space $\mathcal{S}_c(\Psi(\mathcal{A}))$ is Baire.
\end{example}

\proof
Let $D = \N$ and for every  $A \in \mathcal{A}$ we define $X_A = A \cup \{A\}$. Theorem  \ref{SCFCB} asserts that  the space $\mathcal{S}_c(X_A)$ is Baire for each $A \in \mathcal{A}$. Since $D$ and the family $\{X_A : A \in \mathcal{A}\}$ satisfy all the conditions of Lemma \ref{CSCB}, we conclude that  $\mathcal{S}_c(X)$ is Baire.
\endproof

We point out that  the spaces just considered above  have a dense subset of isolated points. The next result shows that without the presence of isolated points in
$X$, the space $\mathcal{S}_c(X)$ can never have the Baire property (in the paper, \cite[Th. 3.2]{sal-yas}, it is shows that if $X$ is a crowded metric space, then $\mathcal{S}_c(X)$ is not Baire).

\begin{theorem}
If $X$ is crowded, then $\mathcal{S}_c(X)$ is meager.
\end{theorem}

\proof As $X$ is crowded, we have that the family of all canonical nonempty open subsets $\left\langle \mathcal{U}\right\rangle$ of $\mathcal{S}_c(X)$ such that $\left|\mathcal{U}\right| \geq 2$ is  a base for $\mathcal{S}_c(X)$. Thus,  in virtue of Lemma \ref{unionopenbaire} $(2)$, it suffices to prove that every such a open set is meager. Indeed,
fix a canonical open set $\left\langle \mathcal{U}\right\rangle$ such that $\left|\mathcal{U}\right| \geq 2$ and  $\left\langle \mathcal{U}\right\rangle = \bigcup\{\mathcal{N}(U,n) : U \in \mathcal{U} \textnormal{ and } n \in \N^+\}$, where $\mathcal{N}(U, n) = \{S \in \left\langle \mathcal{U}\right\rangle : \left|S \setminus U\right| = n\}$ for each  $U \in \mathcal{U}$ and $n \in \N^+$. Let us prove that each set $\mathcal{N}(U,n)$ is nowhere dense. Indeed, pick $U \in \mathcal{U}$ and $n \in \N^+$. Let $O$ be a nonempty open set of $\mathcal{S}_c(X)$. Choose a canonical nonempty open set  $\left\langle \mathcal{V}\right\rangle \subseteq O$ so that  $\left\langle \mathcal{V}\right\rangle \cap \mathcal{N}(U,n) \not= \emptyset$. Since $\left|\mathcal{U}\right| \geq 2$, we can find  $U_0 \in \mathcal{U} \setminus \{U\}$ and $V_0 \in \mathcal{V}$ such that $V_0 \cap U_0 \not= \emptyset$. Since  $X$ is crowded, so we can find a family $\mathcal{W}$ of disjoint nonempty open subsets of $X$ such that $\bigcup\mathcal{W} \subseteq V_0 \cap U_0$ and $\left|\mathcal{W}\right| = n + 1$. Then $\left\langle \mathcal{V} \cup \mathcal{W}\right\rangle \subseteq \left\langle \mathcal{V}\right\rangle \subseteq O$ and $\left\langle \mathcal{V} \cup \mathcal{W}\right\rangle  \cap \mathcal{N}(U,n) = \emptyset$. So,  $\mathcal{N}(U, n)$ is nowhere dense. Therefore, $\left\langle \mathcal{U}\right\rangle$ is meager.
\endproof

The presence of a dense set of isolated points in the above examples is not a causality, as shows the next corollary.

\begin{corollary}\label{SCB}
If $\mathcal{S}_c(X)$ is Baire, then $X$ must have a dense set of isolated points. In particular, $X$ has the Baire property.
\end{corollary}

\proof
Assume that the set of isolated points of $X$ is not dense in $X$. Thus, we can find a nonempty crowded open set $U \subseteq X$. It then follows that $\mathcal{S}_c(U) = \left\langle \mathcal{U}\right\rangle$ is a nonempty  meager open subset of $\mathcal{S}_c(U)$, which is a  contradiction. Therefore,   $X$ must have a dense set of isolated points.
\endproof

As we have seen after Theorem \ref{alpha}, if   $\mathcal{F}_{\mathcal{Q}}$ is the $FAN$-filter, then $\mathcal{S}_c(\xi(\mathcal{F}_{\mathcal{Q}}))$ is meager. Thus, the hyperspace of nontrivial convergence sequence of a Baire space with a dense set of isolated points is not necessarily Baire. To get a space $X$ with $\mathcal{S}_c(X)$ of the second category but not Baire is now easy. Consider the disjoint union  $X = Y \oplus Z$ where $\mathcal{S}_c(Y)$ is Baire and $\mathcal{S}_c(Z)$ is meager. It follows from the equality $\mathcal{S}_c(X) = \mathcal{S}_c(Y) \oplus \mathcal{S}_c(Z) \oplus\left\langle Y,Z\right\rangle$ that $\mathcal{S}_c(X)$ is of the second category but not Baire.

\medskip

Let us remarks  that $\mathcal{S}_c(X)$ can never be countably compact. In fact, if $S \in \mathcal{S}_c(X)$ then $\{S \setminus F : F \in \mathcal{F}(S \setminus\{x_S\})\}$ converges to $\{x_S\}$ in $\mathcal{CL}(X)$, but $\{x_S\} \not \in \mathcal{S}_c(X)$. It follows that $\{S \setminus F : F \in \mathcal{F}(S \setminus\{x_S\})\}$ is closed and discrete in $\mathcal{S}_c(X)$. On the other hand, we have proved that $\mathcal{S}_c(X)$ may have the
Baire property. In this way, is natural to analyze  the pseudocompactness on $\mathcal{S}_c(X)$.

\begin{theorem}\label{SCP}
If $\mathcal{S}_c(X)$ is pseudocompact, then $X$ has a relatively countably compact dense set of isolated points, every finite power of $X$ is pseudocompact, and every $G_\delta$-point in $X$ must be isolated.
\end{theorem}

\proof Assume that $\mathcal{S}_c(X)$ is pseudocompact and let $D$ be the set of isolated points of $X$.
  Hence, $\mathcal{S}_c(X)$ has the Baire property and so, by applying \ref{SCB}, we obtain that $D$ is dense in $X$.
  We claim  that $D$ is relatively countably compact in $X$. Assume that  the contrary, $D$ is not relatively countably compact. Then there exists a countable infinite set $N \subseteq D$ which is clopen $X$. For each $F = \{x_1,\ldots,x_n\} \in \mathcal{F}(N)$ consider the  clopen set  $\mathcal{U}_F = \{X \setminus N,\{x_1\},\ldots,\{x_n\}\}$ and let $\mathcal{U}_\emptyset = \{X \setminus \N\}$. Notice that $\mathcal{S}_c(X) = \bigoplus\{\left\langle \mathcal{U}_F\right\rangle : F \subseteq N \textnormal{ is finite}\}$, but this contradicts the pseudocompactness of $\mathcal{S}_c(X)$. Thus, we have that $D$ is a relatively countably compact subset of $X$.

Now, we shall verify that $X$ is pseudocompact. Suppose that there exists  an infinite  family $\{U_n : n \in \N\}$ of nonempty open subsets of $X$ such that $\textnormal{cl}_X(U_{n+1}) \subseteq U_n$ and $\bigcap\{U_n : n \in \N\} = \emptyset$. Since $D$ is discrete, dense and relatively countably compact in $X$, $X$ is pseudocompact. Besides, since $X$ is Fr\'echet, by Theorem 3.10.26 of \cite{eng} and induction,  $X^n$ is pseudocompact for all $n \in \N$.

Finally assume that $X$ has a non-isolated $G_\delta$-point $x$. By the regularity of $X$, we can find a family $\{U_n : n \in \N\}$ of nonempty open subsets of $X$ such that $\textnormal{cl}(U_{n+1}) \subseteq U_n$ and $\bigcap\{U_n : n \in \N\} = \{x\}$. Then, we must have that  $U_n$ is an infinite sets for all $n \in \mathbb{N}$. So,  $\{\left\langle U_n\right\rangle : n \in \N\}$ is a family of nonempty open subsets of $\mathcal{S}_c(X)$ such that $\textnormal{cl}(\left\langle U_{n+1}\right\rangle) \subseteq \left\langle U_n\right\rangle$ and $\bigcap\{\left\langle U_n \right\rangle : n \in \N\} = \emptyset$, contradicting the pseudocompactness of $\mathcal{S}_c(X)$. Therefore,  the point $x$ is  isolated.
\endproof

There are plenty of spaces satisfying the conclusions from Theorem \ref{SCP}, as a concrete example the Alexandroff duplicate of the $\Sigma$-product of $\omega_1$ copies of the discrete space of cardinality two. However, we do not know  whether this space has hyperspace of convergent sequences pseudocompact. In a more general setting, we have the following question.

\begin{question}\label{q1} Is there a space $X$ for which $\mathcal{S}_c(X)$ is pseudocompact?
\end{question}

We will see in the next theorem that the Alexandroof duplicated will provide examples of spaces with Baire hyperspace of convergent sequences. For the fist example, let us prove a preliminary lemma.

\medskip

For the next lemma, for a space $X$ we let $\A(X) = X \times \{0, 1\}$ denote its Alexandroff duplicate, where $X\times \{1\}$ is discrete. For each $U \subseteq X$,  we define $\widehat{U} = U \times \{0, 1\}$.

\begin{lemma}\label{pibase}
If $\mathcal{B}$ is a $\pi$-base of a space  $X$ and $\mathcal{B}^\ast$ is a $\pi$-base of the set of all non-isolated points of $X$ consisting of non-discrete sets, then the family of all canonical sets of the form
$$\langle \{B\} \cup \mathcal{U} \rangle,$$ where $B \in \mathcal{B}^\ast$ and $\mathcal{U} \subseteq \mathcal{B}$, is a $\pi$-base of $\mathcal{S}_c(\A(X))$ consisting of nonempty open sets.
\end{lemma}

\proof Let $\left\langle \mathcal{V}\right\rangle$ be a canonical nonempty open subset of $\mathcal{S}_c(\A(X))$. Note that there exists a non discrete set $\mathcal{V}_0 \in \mathcal{V}$. Pick $B \in \mathcal{B}^\ast$ such that $B \subseteq V_0$ and $B_V \in \mathcal{B}$ such that $B_V \subseteq V$ for each $V \in \mathcal{V} \setminus \{V_0\}$. Then, it is clear that $\emptyset \neq  \langle \{B\} \cup \mathcal{U} \rangle \subseteq \left\langle \mathcal{V}\right\rangle$.
\endproof

For the next result, the diameter of a subset $A$ of a metric space $(X,d)$ will be denoted by  $\delta(A) := \sup \{ d(x,y) : x, y \in A \}$.

\begin{theorem}\label{comme}
If $X$ is a complete metrizable crowded space, then $\mathcal{S}_c(\A(X))$ has the Baire property.
\end{theorem}

\proof Equipped $X$ with a complete compatible metric. Set $\mathcal{B} = \{\{(x,1)\} : x \in X\}$ and $\mathcal{B}^\ast = \{ \widehat{U} : U \subseteq X  \textnormal{ is open} \}$. Note that $\mathcal{B}$ is a $\pi$-base for $\mathcal{A}(X)$, because of $X$ does not have isolated points, and $\mathcal{B}^\ast$ is a base for $X \times \{0\}$. Suppose that  $\{\mathcal{D}_i : i\in \N \}$ is a decreasing sequence of dense open subsets of $\mathcal{S}_c(\A(X))$ and let  $\left\langle \mathcal{U}\right\rangle$ be a canonical open subset of $\mathcal{S}_c(\A(X))$, where $\mathcal{U}  = \{U_1,\ldots,U_k\}$  and $\{U_1,\ldots,U_k\}$ are pairwise disjoint  nonempty open subsets of $\A(X)$. By inductively applying Lemma \ref{pibase}, we can find
   a strictly increasing sequence $(n_i)_{i \in \mathbb{N}}$ in $\mathbb{N}$, a sequence of points $(x_i)_{i \in \mathbb{N}}$, and for every $i \in \mathbb{N}$,  a nonempty open subset $W_i$ of $X$ such that:
\begin{enumerate}[$(a)$]
\item $cl_X(W_{i+1}) \subseteq W_i$, for every $i \in \N$.

\item $\delta(W_i) < \frac{1}{2^i}$, for every $i \in \N$.

\item $\widehat{W}_i \cap \{(x_1,1),\dots, (x_{n_i},1)\} = \emptyset$, for every $i \in \N$.

\item $\left\langle \widehat{W}_i,\{(x_1,1)\},\dots,\{(x_{n_i},1)\} \right\rangle \subseteq \left\langle \mathcal{U}\right\rangle \cap \mathcal{D}_i$, for every $i \in \N$.

\item $x_j \in W_i$, for all $i, j \in \N$ with $n_i < j$.
\end{enumerate}
It follows from $(b)$ and $(e)$ that  $(x_i)_{i \in \mathbb{N}}$ is a Cauchy sequence and since $X$ is complete, there is $x \in \bigcap_{i \in \N}W_i$ such that $x_i \to x$. Consider the sequence $s = \{x\} \cup \{ x_i : i \in \N\}$. It is evident from the construction that
$s \in  \left\langle \mathcal{U}\right\rangle \cap \big( \bigcap_{i \in \N}\mathcal{D}_i \big)$. Therefore,  $\mathcal{S}_c(\A(X))$ has the Baire property.
\endproof

Theorem \ref{comme} can be generalized as follows.

\begin{theorem}\label{gecomme}
If $Z$ is a $\Sigma$-product of completely metrizable spaces and $Z$ is crowded, then $\mathcal{S}_c(\mathcal{A}(Z))$ has the Baire property.
\end{theorem}

\proof Assume that $Z = \{x \in X : \left| \textnormal{suppt}(x)\right| \leq \omega \}$, where $X = \prod_{i \in I}X_i$ is a product of completely metrizable spaces, $a = (a_i)_{i \in I} \in X$  is a fixed point and $\textnormal{suppt}(x) := \{ i \in I : x_i \neq a_i \}$ for each $x \in X$.
For each $i \in I$ we equipped $X_i$ with a a complete metric. Set $\mathcal{B} = \{\{(x,1)\} : x \in Z\}$
and note that $\mathcal{B}$ is a $\pi$-base for $\mathcal{A}(Z)$ because of $Z$ does not have isolated points.
Now, for each $n \in \N$, consider  the family $\mathcal{B}_n$ of all nonempty canonical open subsets $B$ of $Z$ such that the projection  $\pi_i[B]$ has diameter smaller than $\frac{1}{2^n}$, for all $i \in \textnormal{supp}(B)\footnote{For a canonical open set $B$ of $X$, we define $\textnormal{supp}(B) = \{ i \in I : \pi_i[B] \neq X_i\}$}$. It is clear that each $\mathcal{B}_n$ is a base for the space $Z$ consisting of crowded sets. As a consequence $\widehat{\mathcal{B}}_n = \{ \widehat{B} : B \in \mathcal{B}_n\}$ is a family of open non-discrete sets which is a $\pi$-base at each non-isolated point of $\mathcal{A}(Z)$ for each $n \in \N$. To show that $\mathcal{S}_c(\mathcal{A}(Z))$ is a Baire space suppose that $\{\mathcal{D}_n : n \in \N \}$ is a family of open dense subsets of $\mathcal{S}_c(\mathcal{A}(Z))$. Fix an arbitrary nonempty canonical open subset $\left\langle \mathcal{U}\right\rangle$ of $\mathcal{S}_c(\mathcal{A}(Z))$. Let $\{N_n : n \in \N\}$ be a partition of $\N$ in infinite subsets. We will construct recursively, for each $n \in \N$, a set $B_n \in \mathcal{B}_n$ and a finite subset $\mathcal{U}_n$ of $\mathcal{B}$ as follows:

\smallskip

 By using the Lemma \ref{pibase}, we can find a set $B_1 \in \mathcal{B}_1$ and a finite subset $\mathcal{U}_1$ of $\mathcal{B}$ such that $\langle\{ \widehat{B}_1\}\cup \mathcal{U}_1\rangle$ is a canonical open set contained in $\left\langle \mathcal{U}\right\rangle \cap \mathcal{D}_1$. We can may assume that the cardinality of $\mathcal{U}_1$ is at least two. Let $F_1 = \bigcup\mathcal{U}_1$ and $A_1 = \bigcup\{\textnormal{suppt}(x) : (x,1) \in F_1 \}$.  Enumerate $A_1$ as $\{i_m : m \in N_1\}$, we may repeat the elements of $F_1$ if it is necessary. Assume that for each $k \leq n$ we have defined $B_k \in \mathcal{B}_k$ and $\mathcal{U}_k \in \mathcal{B}$ satisfying the corresponding conditions $(1)-(4)$ below. Before to procedure to the next induction step, for every $k \leq n - 1$, we define $F_{k + 1} = \bigcup(\mathcal{U}_{k + 1} \setminus \mathcal{U}_k)$ and enumerate $A_k := \bigcup\{\textnormal{suppt}(x) : (x,1) \in F_k \}$ as  $\{i_{m} : m \in N_{k}\}$ allowing repetition if it necessary. Them, by applying Lemma \ref{pibase} again, we can find a set $B_{n + 1} \in \mathcal{B}_{n+1}$ and a finite subset $\mathcal{U}_{n + 1}$ of $\mathcal{B}$ so that:
\begin{enumerate}
\item $\langle \{\widehat{B}_{n+1}\}\cup \mathcal{U}_{n+1}\rangle$ is a canonical open set contained in $\langle \{\widehat{B}_{n}\}\cup \mathcal{U}_{n}\rangle \cap \mathcal{D}_{n+1}$,

\item $\textnormal{cl}_Z(B_{n+1}) \subseteq B_{n}$,

\item $\mathcal{U}_{n+1} \setminus \mathcal{U}_{n}$ has at least two elements, and

\item $\{i_m : m \in \big(\bigcup_{k \leq n}N_k\big) \cap n  \} \subseteq \textnormal{supp}(B_{n+1})$.
\end{enumerate}
Thus, we have defined $B_n$ and $\mathcal{U}_n$ for each $n \in \N$. It follows from $(4)$ that $A := \bigcup_{n \in \N}A_n = \{i_m : n \in \N\} \subseteq \bigcup_{n \in \N}\textnormal{supp}(B_{n}) =: C$. On the other hand, we deduce from $(1)$ that if $G_n := \{x : (x,1) \in F_n\}$ for every $n \in \N$, then $G_{n+1} \subseteq B_n$.  Fix $i \in C$. Since $\{\textnormal{cl}(\pi_i[B_n]) : n \in \N\}$ is a deceasing sequence of closed subsets of $X_i$ whose  diameters converge to $0$,  there exists a unique point $b_i \in \bigcap_{n \in \N} \textnormal{cl}[\pi_i(B_n)]$.  Next we proceed to define $z \in Z$ as follows:
$$
z_i =
\begin{cases}
b_i &  \textnormal{if } i \in C\\
a_i   &  \textnormal{if } i \in I \setminus C.
\end{cases}
$$
Our desired sequence is $S = \{(z,0)\} \cup \big(\bigcup_{n \in \N}F_n\big)$. Indeed, it is evident from $(3)$ that $\bigcup_{n \in \N}F_n$ is infinite and discrete. To see that the sequence $\bigcup_{n \in \N}F_n$ converges to $(z,0)$. First, when $i\in I \setminus C$  the fact that $\bigcup\{\textnormal{suppt}(x) : x \in \bigcup_{n \in \N}G_{n} \} = A \subseteq C$ implies that
 $\pi_i[\bigcup_{n \in \N}G_{n}] = \{a_i\}$ and so $\pi_i[\bigcup_{n \in \N} G_{n}]$ trivially converges to $a_i = z_i$. Secondly, when $i \in C$ it is clear from the construction that $\pi_i[\bigcup_{n \in \N} G_{n}]$ converges to $b_i = z_i$. It then follows that $\bigcup_{n \in \N} G_{n}$ converges to $z$. Since $Z$ is crowded, we conclude that $\bigcup_{n \in \N}F_n$ converges to $(z,0)$.
\endproof

It is evident that Theorem \ref{comme} follows directly from Theorem \ref{gecomme}, but we decide to include the proofs of booth since the proof of the former  is illustrative and didactic, and after reading it  the proof of the later will be more understandable.

\medskip

We list several open questions which are closed related to the results of this article.

\begin{question} Is there a space $X$ for which $\mathcal{S}_c(\A(X))$ is pseudocompact?
\end{question}

Examples \ref{ex1} and \ref{ex2} suggest the following question.

\begin{question} Characterize the Baire spaces $X$'s for which  $\mathcal{S}_c(X)$ is Baire?
\end{question}

Based on Corollary \ref{SCB} we have formulated the next question.

\begin{question} What are the properties of a space $X$ when   $\mathcal{S}_c(X)$ is second category?
\end{question}

Following the paper \cite{gg}, we say that a space $X$ is {\it weakly pseudocompact} if it is $G_\delta$-dense in some compactification. We know that every pseudocompact space is weakly pseudocompact  and every weakly pseudocompact space is Baire (for the details of these facts see \cite{gg}).  Thus, we may weaken Question \ref{q1} as:

\begin{question}\label{q1} Is there a space $X$ for which $\mathcal{S}_c(X)$ is weakly pseudocompact?
\end{question}


\end{document}